\input amstex
\input amsppt.sty
\magnification=\magstep1
\hsize=30truecc
\vsize=22.2truecm
\baselineskip=16truept
\TagsOnRight
\nologo
\pageno=1
\topmatter
\def\N{\Bbb N}
\def\Z{\Bbb Z}

\def\l{\left}
\def\r{\right}
\def\b{\bigg}

\def\({\b(}
\def\[{\b[}
\def\){\b)}
\def\]{\b]}

\def\t{\text}
\def\f{\frac}
\def\mo{\roman{mod}}
\def\ord{\roman{ord}}

\def\bi{\binom}
\def\eq{\equiv}

\def\ls{\leqslant}
\def\gs{\geqslant}

\def\da{\delta}

\def\mm#1#2#3{\thickfracwithdelims\{\}\thickness0{#1}{#2}_{#3}}
\def\Proof{\noindent{\it Proof}}
\def\Remark{\medskip\noindent{\it Remark}}

\hbox {Int. J. Number Theory 4(2008), no.\,2, 155--170.}
\bigskip
\title Lucas-type congruences for cyclotomic $\psi$-coefficients\endtitle
\author Zhi-Wei Sun$^1$ and Daqing Wan$^2$\endauthor
\leftheadtext{Zhi-Wei Sun and Daqing Wan}
\affil $^1$Department of Mathematics, Nanjing University
\\ Nanjing 210093, People's Republic of China
\\zwsun\@nju.edu.cn
\\ {\tt http://math.nju.edu.cn/$\sim$zwsun}
\medskip
$^2$Department of Mathematics, University of California\\Irvine, CA 92697-3875, USA
\\dwan\@math.uci.edu
\\ {\tt http://math.uci.edu/$\sim$dwan}
\endaffil
\thanks
\indent The
first author is supported by the National Science Fund for
Distinguished Young Scholars in China (grant no. 10425103). The
second author is partially supported by NSF.
\endthanks
\abstract
Let $p$ be any prime and $a$ be a positive integer.
For $l,n\in\{0,1,\ldots\}$ and $r\in\Z$,
the normalized cyclotomic $\psi$-coefficient
$$\mm nr{l,p^a}:=p^{-\l\lfloor\f{n-p^{a-1}-lp^{a}}{p^{a-1}(p-1)}\r\rfloor}
\sum_{k\eq r\,(\mo\ p^{a})}(-1)^k\bi nk\bi{\f{k-r}{p^a}}l$$
is known to be an integer. In this paper, we show that
this coefficient behaves like binomial coefficients and satisfies
some Lucas-type congruences. This implies that a
congruence of Wan is often optimal, and two conjectures of Sun and Davis
are true.
\endabstract
\endtopmatter

\document

\heading{1. Introduction}\endheading

As usual, the binomial coefficient $\bi x0$ is regarded as $1$.
For $k\in\Z^+=\{1,2,\ldots\}$, we define
$$\bi xk=\f{x(x-1)\cdots(x-k+1)}{k!}$$
and adopt the convention $\bi x{-k}=0$.

The following remarkable result was established
by A. Fleck (cf. [D, p.\,274]) in the case $l=0$ and $a=1$,
by C. S. Weisman [We] in the case $l=0$, and by D. Wan [W] in the general case
motivated by his study of the $\psi$-operator related to Iwasawa theory.

\proclaim{Theorem 1.0} Let $p$ be a prime and $a\in\Z^+$.
Then, for any $l,n\in\N=\{0,1,2,\ldots\}$ and $r\in\Z$, we have
$$C_{l, p^a}(n, r):=
\sum_{k\eq r\,(\mo\ p^{a})}(-1)^k\bi nk\bi{(k-r)/p^{a}}l\in
p^{\lfloor\f{n-p^{a-1}-lp^{a}}{\phi(p^{a})}\rfloor}\Z,$$
where $\phi$ is Euler's totient function and $\lfloor\cdot\rfloor$
is the greatest integer function.
\endproclaim

The above integers $C_{l, p^a}(n, r)\ (l=0,1,\ldots)$ arise naturally as the
coefficients of the $\psi$-operator acting on the cyclotomic
$\varphi$-module. We briefly review this connection.
Let $A={\Z}_p [\![T]\!]$
be the formal power series ring over the ring
of $p$-adic integers.
The $\Z_p$-linear Frobenius map $\varphi$ acts on the ring $A$ by
$$\varphi(T) =(1+T)^p-1.$$
Equivalently, $\varphi(1+T)
=(1+T)^p$. This map $\varphi$ is injective and of degree $p$. This
implies that
$\{1, T,\ldots,T^{p-1}\}$ and
 $\{1, 1+T,\ldots, (1+T)^{p-1}\}$ are bases of
$A$ over the subring $\varphi(A)$.
The operator $\psi:A\to A$ is defined by
$$\psi(x) = \psi\(\sum_{i=0}^{p-1}(1+T)^i\varphi(x_i)\)
=x_0=\frac{1}{p} \varphi^{-1}(\text{Tr}_{A/\varphi(A)}(x)),$$
where $x: A \to A$ denotes the multiplication by $x$ as a
$\varphi(A)$-linear map. Note that $\psi$ is a one-sided inverse of $\varphi$,
namely $\psi\circ \varphi=I\not=\varphi\circ\psi$.
The pair $(A, \varphi)$ is the
cyclotomic $\varphi$-module. The $\psi$-operator plays a basic role in
$L$-functions of $F$-crystals,
Fontaine's theory of $(\varphi, \Gamma)$-modules, Iwasawa theory,
$p$-adic $L$-functions and $p$-adic Langlands correspondence.

For a positive integer $a$, let $\psi^a$ be the $a$-th iteration
of $\psi$ acting on the ring $A$.
As mentioned in [W, Lemma 4.2],
it is easy to check
that for any $n\in\N$ and $r\in\Z$ we have
$$\psi^a \l(\f{T^n}{(1+T)^{r}}\r) = (-1)^n \sum_{l=0}^{\infty} T^l C_{l, p^a}(n, r).$$
To understand the $\psi^a$-action, it is thus essential to
understand the $p$-adic property of the cyclotomic $\psi$-coefficients
$C_{l, p^a}(n,r)\ (l=0,1,\ldots)$. This was the main motivation in [W], where
the congruence in Theorem 1.0 was proved. Note that a somewhat
weaker estimate for the cyclotomic $\psi$-coefficient $C_{l, p}(n,0)$
was independently given by Colmez [C, Lemma 1.7] in his work
on $p$-adic Langlands correspondence. The cyclotomic $\psi$-coefficient
also arises from computing the homotopy $p$-exponent of the
special unitary group $\t{SU}(n)$ (cf. [DS]).

To understand how sharp the congruence in Theorem 1.0 is,
we define the {\it normalized cyclotomic $\psi$-coefficient}
$$\mm nr{l,p^{a}}:=
p^{-\l\lfloor\f{n-p^{a-1}-lp^{a}}{\phi(p^{a})}\r\rfloor}
\sum_{k\eq r\,(\mo\ p^{a})}(-1)^k\bi nk\bi{(k-r)/p^{a}}l.\tag1.0$$
Surprisingly it has many properties similar to
properties of the usual
binomial coefficients.

The classical Lucas theorem states that if $p$ is a prime and
$n,r,s,t$ are nonnegative integers with $s,t<p$ then
$$\bi{pn+s}{pr+t}\eq\bi nr\bi st\ \ (\mo\ p).$$
It can also be interpreted as a result about cellular automata (cf. [Gr]).
There are various extensions of this fundamental theorem,
see, e.g., [DW], [HS], [P] and [SD].
Our first result is the following new analogue of Lucas' theorem.

\proclaim{Theorem 1.1}
Let $p$ be any prime, and let $r\in\Z$ and $a,l,n,s,t\in\N$ with
$a\gs2$ and $s,t<p$. Then we have the congruence
$$\mm{pn+s}{pr+t}{l,p^{a+1}}\eq(-1)^t\bi st\mm{n}{r}{l,p^{a}}\ \ (\mo\ p);\tag1.1$$
in other words,
$$\align &p^{-\l\lfloor\f{n-p^{a-1}-lp^{a}}{\phi(p^{a})}\r\rfloor}
\sum_{k\eq r\,(\mo\ p^{a})}(-1)^{pk}\bi{pn+s}{pk+t}\bi{(k-r)/p^{a}}l
\\\eq&p^{-\l\lfloor\f{n-p^{a-1}-lp^{a}}{\phi(p^{a})}\r\rfloor}
\sum_{k\eq r\,(\mo\ p^{a})}(-1)^k\bi nk\bi st\bi{(k-r)/p^{a}}l\ (\mo\ p).
\endalign$$
\endproclaim

\Remark\ 1.1. Theorem 1.1 in the case $l=0$ is equivalent to
Theorem 1.7 of Z. W. Sun and D. M. Davis [SD].
Under the same conditions of Theorem 1.1, Sun and Davis [SD]
established another congruence of Lucas' type:
$$\align &\f1{\lfloor n/p^{a-1}\rfloor!}
\sum_{k\eq r\,(\mo\ p^{a})}(-1)^{pk}\bi{pn+s}{pk+t}\l(\f{k-r}{p^{a-1}}\r)^l
\\\eq&\f1{\lfloor n/p^{a-1}\rfloor!}
\sum_{k\eq r\,(\mo\ p^{a})}(-1)^k\bi nk\bi st\l(\f{k-r}{p^{a-1}}\r)^l\ (\mo\ p).
\endalign$$
\medskip

Note that $a\gs 2$ is assumed in Theorem 1.1. To get a complete
result, we need to handle the case $a=1$ as well, which is more
subtle. In fact, concerning the exceptional case $a=1$, Sun and
Davis [SD] made the following conjecture (for $l=0$). Note also
that [S02] contains a closed formula for $\mm nr{0,2^2}$ with
$n\in\N$ and $r\in\Z$.

\proclaim{Conjecture {\rm ([SD, Conjecture 1.2])}} Let $p$ be any prime,
and let $n\in\N$, $r\in\Z$ and $s\in\{0,\ldots,p-1\}$.
If $p\mid n$ or $p-1\nmid n-1$, then
$$\mm{pn+s}{pr+t}{0,p^2}\eq(-1)^t\bi st
\mm nr{0,p}\ \ (\mo\ p)$$
for every $t=0,\ldots,p-1$.
When $p\nmid n$ and $p-1\mid n-1$, the least nonnegative residue of
$\mm{pn+s}{pr+t}{0,p^2}$ modulo $p$
does not depend on $r$ for each integer $t\in(s,p-1]$, moreover these residues
form a permutation of $1,\ldots,p-1$ if $s=0$ and $n\not=1$.
\endproclaim

We get the following general result
for $a=1$ and all $l\in\N$ from which the above conjecture follows.

\proclaim{Theorem 1.2} Let $p$ be a prime, $l,n\in\N$, $r\in\Z$
and $s,t\in\{0,\ldots,p-1\}$. If $p\mid n$, or $p-1\nmid n-l-1$,
or $s=p-1$, or $s=2t$ and $p\not=2$, then
$$\mm{pn+s}{pr+t}{l,p^2}\eq(-1)^t\bi st\mm nr{l,p}\ \ (\mo\ p).\tag1.2$$
When $p\nmid n$, $p-1\mid n-l-1$ and $t\in(s,p-1]$, we have
$$\mm{pn+s}{pr+t}{l,p^2}
\eq\cases(-1)^{s+\f{n-l-1}{p-1}}\,\f nt\bi{\f{n-l-1}{p-1}-1}l/\bi{t-1}s\  (\mo\ p)
&\t{if}\ n>l+1,\\0\ (\mo\ p)&\t{if}\ n\ls l+1.\endcases\tag1.3$$
\endproclaim

From Theorem 1.2 we can also deduce the following result conjectured by Sun and Davis
(cf. [SD, Remark 1.4]) as a complement to Theorem 1.5 of [SD].

\proclaim{Corollary 1.3} Let $p$ be any prime, and let $l,n\in\N$ and $r\in\Z$.
Then
$$T^{(p)}_{l,2}(n,r)\eq(-1)^{\{r\}_p}\bi{\{n\}_p}{\{r\}_p}T^{(p)}_{l,1}
\l(\l\lfloor\f np\r\rfloor,\,\l\lfloor\f rp\r\rfloor\r)\ \ (\mo\ p),\tag1.4$$
where
$$T^{(p)}_{l,a}(n,r):=\f{l!p^l}{\lfloor n/p^{a-1}\rfloor!}
\sum_{k\eq r\,(\mo\ p^a)}(-1)^k\bi nk\bi{(k-r)/p^a}l\ \t{for}\ a\in\Z^+,$$
and we use $\{x\}_m$ to denote the least
nonnegative residue of an integer $x$ modulo $m\in\Z^+$.
\endproclaim

When $s=t=0$, the Lucas-type congruences in Theorems 1.1 and 1.2
can be further improved unless $p=2$ and $2\nmid n$.
Namely, we have the following result.

\proclaim{Theorem 1.4} Let $p$ be a prime, and let $a,n\in\Z^+$, $l\in\N$ and $r\in\Z$. Then
$$\aligned\ord_p\(\mm{pn}{pr}{l,p^{a+1}}-\mm nr{l,p^a}\)
\gs\f{p-1}p(2\ord_p(n)+\da),
\endaligned\tag1.5$$
where $\ord_p(n)=\sup\{m\in\N:\,p^m\mid n\}$ and
$$\da=\cases0&\t{if}\ p=2,\\1&\t{if}\ p=3,\\2&\t{if}\ p\gs 5.\endcases$$
\endproclaim

\Remark\ 1.2. Let $p$ be a prime, $a,n\in\Z^+$ and $r\in\Z$.
Substituting $p^{a-1}n$ for $n$ in (1.5), we obtain that
$$\aligned&\ord_p\(\mm{p^a n}{pr}{l,p^{a+1}}-\mm {p^{a-1}n}r{l,p^a}\)
\\\gs&\f{p-1}p(2\ord_p(p^{a-1}n)+\da)\gs\f{p-1}p(2(a-1)+\da).
\endaligned$$
On the other hand, in the case $l=0$
Sun and Davis [SD, Theorem 3.1] proved the congruence
$$\mm{p^a n}{pr}{0,p^{a+1}}\eq\mm{p^{a-1}n}r{0,p^a}\ \l(\mo\ p^{(2-\da_{p,2})(a-1)}\r)$$
(where the Kronecker symbol $\da_{i,j}$ takes $1$ or $0$ according as $i=j$ or not)
and they conjectured that the exponent $(2-\da_{p,2})(a-1)$
can be replaced by $2a-\da_{p,3}=2(a-1)+\da$ when $p\not=2$.
\medskip

Here is one more result, which shows that Theorem 1.0 is often sharp.

\proclaim{Theorem 1.5} Let $p$ be any prime, and let $a\in\Z^+$ and $l\in\N$.
If $n=(l+1)p^{a-1}-1+m\phi(p^a)$ for some $m\in\Z^+$, then
$$\mm nr{l,p^a}\eq(-1)^{m-1}\bi{m-1}l\ \ (\mo\ p)\quad\t{for all}\ r\in\Z.\tag1.6$$
\endproclaim

\Remark\ 1.3. Theorem 1.5 in the case $l=0$ was first obtained by
Weisman [We] in 1977. Given $l\in\Z^+$,
for any integer $m>l$ with $m\eq l+1\ (\mo\ p^{\lfloor\log_p l\rfloor+1})$
we have $\bi{m-1}l\eq\bi ll=1\ (\mo\ p)$ by Lucas' theorem.
\medskip

In the next section we include a new
proof of Theorem 1.0 of a combinatorial nature.
In Section 3 we will show Theorem 1.1.
Theorems 1.2 and Corollary 1.3 will be proved in Section 4.
Section 5 is devoted to proofs of Theorems 1.4 and 1.5.
Instead of the $\psi$-operator, we use combinatorial arguments
throughout this paper.

\heading{2. A combinatorial proof of Theorem 1.0}\endheading

\proclaim{Lemma 2.1} Let $a,b\in\Z$ and $m\in\Z^+$. Then
$$\l\lfloor\f am\r\rfloor+\l\lfloor\f bm\r\rfloor+1
-\l\lfloor\f{a+b+1}m\r\rfloor\in\{0,1\}.\tag2.1$$
\endproclaim
\Proof. Observe that
$$\l\lfloor\f{a+b+1}m\r\rfloor=\l\lfloor\f am\r\rfloor+\l\lfloor\f bm\r\rfloor
+\l\lfloor\f{\{a\}_m+\{b\}_m+1}m\r\rfloor. $$ The last term is
obviously either $0$ or $1$, so (2.1) follows. \qed

\proclaim{Lemma 2.2} Let $l,m,n\in\Z^+$ and $r\in\Z$. Then we have
$$\align&\sum_{k\eq r\,(\mo\ m)}(-1)^k\bi nk\bi{(k-r)/m}l
-\bi{\lfloor(n-r)/m\rfloor}l\sum_{m\mid k-r}(-1)^k\bi nk
\\&=-\sum_{j=0}^{n-1}\bi nj\sum_{m\mid i-r}(-1)^i\bi ji
\sum_{m\mid k-r_j}(-1)^k\bi{n-j-1}{k}
\bi{(k-r_j)/m}{l-1},
\endalign$$
where $r_j=r-j+m-1$.
\endproclaim
\Proof. Note that $\bi{x+1}l-\bi{x}l=\bi{x}{l-1}$.
So Lemma 2.2 is just Lemma 3.3 of [DS] in the case $f(x)=\bi xl$. \qed

\bigskip
\noindent{\it Proof of Theorem 1.0}. We use induction on $l+n$.

The case $n=0$ is trivial.
The case $l=0$ was handled by Weisman [W] (see also [S06]).

Now let $l$ and $n$ be positive, and assume that $\mm{n'}{r'}{l',p^a}\in\Z$
whenever $l',n'\in\N$, $l'+n'<l+n$ and $r'\in\Z$.
By Lemma 2.2,
$$\align&\mm nr{l,p^a}-\bi{\lfloor\f{n-r}{p^{a}}\rfloor}l
p^{\l\lfloor \f{n-p^{a-1}}{\phi(p^{a})}\r\rfloor
-\l\lfloor\f{n-p^{a-1}-lp^{a}}{\phi(p^{a})}\r\rfloor}\mm nr{0,p^a}
\\&\ \ =-\sum_{j=0}^{n-1}\bi nj
p^{c_j}\mm jr{0,p^{a}}\mm{n-j-1}{r_j}{l-1,p^a}
\endalign$$
where
$$\align c_j=&\l\lfloor\f{j-p^{a-1}}{\phi(p^{a})}\r\rfloor
+\l\lfloor\f{n-j-1-p^{a-1}-(l-1)p^{a}}{\phi(p^{a})}\r\rfloor
-\l\lfloor\f{n-p^{a-1}-lp^{a}}{\phi(p^{a})}\r\rfloor
\\=&\l\lfloor\f {a_j}{\phi(p^{a})}\r\rfloor
+\l\lfloor\f{b_j}{\phi(p^{a})}\r\rfloor+1
-\l\lfloor\f{a_j+b_j+1}{\phi(p^{a})}\r\rfloor\gs0\quad (\t{by Lemma 2.1})
\endalign$$
with $a_j=j-p^{a-1}$ and $b_j=n-j-1-lp^{a}$.
For any $j=0,1,\ldots,n-1$, both $\mm jr{0,p^a}$ and $\mm{n-j-1}{r_j}{l-1,p^a}$
are integers by the induction hypothesis. Therefore $\mm nr{l,p^a}\in\Z$
by the above.

The induction proof of Theorem 1.0 is now complete. \qed

\Remark\ 2.1. Our proof of Theorem 1.0 can be refined to
show the following recurrence with respect to $l$:
If $p$ is a prime, $a,l,n\in\Z^+$ and $r\in\Z$, then
$$\mm nr{l,p^{a}}\eq-\sum_{j\in J}\bi nj\mm jr{0,p^a}\mm {n-j-1}{r-j+p^a-1}{l-1,p^a}
\ (\mo\ p),$$
where
$$J=\l\{0\ls j\ls n-1:\ \{j-p^{a-1}\}_{\phi(p^{a})}\gs
\{n-(l+1)p^{a-1}\}_{\phi(p^{a})}\r\}.$$

\medskip

\heading{3. Proof of Theorem 1.1}\endheading

We can deduce Theorem 1.1 by using Remark 2.1 along with Theorem 1.7 of [SD].
However, we will present a self-contained proof by a new approach.

\proclaim{Lemma 3.1} Let $d,q\in\Z^+$, $n\in\N$,
$r,t\in\Z$ and $t<d$. Then
$$\aligned&\sum_{j\in\N}(-1)^j
\(\sum_{d\mid k-t}(-1)^k\bi nk\bi{(k-t)/d}j\)
\(\sum_{q\mid i-r}(-1)^{i}\bi ji\bi{(i-r)/q}l\)
\\&\quad=\sum_{k\eq dr+t\,(\mo\ dq)}(-1)^k\bi nk\bi{(k-dr-t)/(dq)}l.
\endaligned\tag3.1$$
\endproclaim

\Proof. Since $t<d$, we have $(k-t)/d\in\N$ for those
$k\in\{0,\ldots,n\}$ with $k\eq t\ (\mo\ d)$. Let $S$ denote the
left-hand side of (3.1). Then
$$S=\sum_{k\eq t\,(\mo\ d)}(-1)^k\bi nk
\sum_{q\mid i-r}\bi{(i-r)/q}l\sum_{j\gs i}(-1)^{j-i}\bi{(k-t)/d}j\bi ji.$$
The inner-most sum has a well-known evaluation (see, e.g., [G, (3.47)] or [GKP, (5.24)]);
in fact, it coincides with
$$\bi{(k-t)/d}i\sum_{j\gs i}(-1)^{j-i}\bi {(k-t)/d-i}{j-i}=\da_{i,(k-t)/d}.$$
Therefore
$$\align S=&\sum_{k\eq t\,(\mo\ d)}(-1)^k\bi nk\sum_{q\mid i-r}
\bi{(i-r)/q}l\da_{i,(k-t)/d}
\\=&\sum_{k\eq dr+t\,(\mo\ dq)}(-1)^k\bi nk\bi{((k-t)/d-r)/q}l
\\=&\sum_{k\eq dr+t\,(\mo\ dq)}(-1)^k\bi nk\bi{(k-dr-t)/(dq)}l.
\endalign$$
This concludes the proof. \qed

\proclaim{Lemma 3.2} Let $p$ be a prime, and let $a\in\Z^+$ and
$l,n\in\N$. Let $r\in\Z$ and $s,t\in\{0,1,\ldots,p-1\}$. If $n=0$
or $s=p-1$ or $\phi(p^a)\nmid n-(l+1)p^{a-1}$, then $(1.1)$ holds;
otherwise,
$$\aligned&\mm{pn+s}{pr+t}{l,p^{a+1}}-(-1)^t\bi st\mm nr{l,p^a}
\\\eq&(-1)^{n-1}\mm{n-1}r{l,p^a}
\mm{pn+s}t{n-1,p}\ (\mo\ p).
\endaligned\tag3.2$$
\endproclaim
\Proof. Applying Lemma 3.1 with $d=p$ and $q=p^a$, we find that
$$\align&\sum_{j\in\N}(-1)^jp^{\l\lfloor\f{pn+s-1-jp}{\phi(p)}\r\rfloor}\mm {pn+s}t{j,p}
p^{\l\lfloor \f{j-p^{a-1}-lp^a}{\phi(p^a)}\r\rfloor}\mm jr{l,p^a}
\\&=p^{\l\lfloor \f{pn+s-p^{a}-lp^{a+1}}{\phi(p^{a+1})}\r\rfloor}\mm {pn+s}{pr+t}{l,p^{a+1}}.
\endalign$$
Thus
$$\mm {pn+s}{pr+t}{l,p^{a+1}}
=\sum_{0\ls j\ls\lfloor\f{pn+s}p\rfloor=n}(-1)^jp^{a_j}\mm jr{l,p^a}\mm {pn+s}t{j,p},$$
where
$$\align a_j=&\l\lfloor\f{pn+s-1-jp}{\phi(p)}\r\rfloor
+\l\lfloor\f{j-p^{a-1}-lp^a}{\phi(p^a)}\r\rfloor
-\l\lfloor\f{pn+s-p^a-lp^{a+1}}{\phi(p^{a+1})}\r\rfloor
\\=&\l\lfloor\f{p(n-j)+s-1}{\phi(p)}\r\rfloor
+\l\lfloor\f{j-p^{a-1}-lp^a}{\phi(p^a)}\r\rfloor
-\l\lfloor\f{n-p^{a-1}-lp^a}{\phi(p^a)}\r\rfloor.
\endalign$$

Observe that
$$\align p^{a_n}\mm{pn+s}t{n,p}=&\sum_{k\eq t\,(\mo\ p)}(-1)^k\bi{pn+s}k\bi{(k-t)/p}n
\\=&(-1)^{pn+t}\bi{pn+s}{pn+t}\bi{(pn+t-t)/p}n
\\\eq&(-1)^{n+t}\bi st\ \ \ (\mo\ p)
\endalign$$
where we have applied Lucas' theorem in the last step.

When $n$ is positive, clearly
$$\align a_{n-1}-\l\lfloor\f s{p-1}\r\rfloor=&1+\l\lfloor\f{n-1-p^{a-1}-lp^a}{\phi(p^a)}\r\rfloor
-\l\lfloor\f{n-p^{a-1}-lp^a}{\phi(p^a)}\r\rfloor
\\=&\cases1&\t{if}\ n\not\eq p^{a-1}+lp^a\eq(l+1)p^{a-1}\ (\mo\ \phi(p^a)),
\\0&\t{otherwise}.\endcases
\endalign$$

Let $0\ls j\ls n-2$. We will see that $a_j\gs n-j-1\gs1$. Since
$$p^a(n-j)+p^{a-1}(s-1)-(n-j)\phi(p^a)=p^{a-1}(n-j+s-1)\gs n-j-1,$$
we have
$$\l\lfloor\f{p(n-j)+s-1}{\phi(p)}\r\rfloor
=\l\lfloor\f{p^a(n-j)+p^{a-1}(s-1)}{\phi(p^a)}\r\rfloor
\gs\l\lfloor\f{n-j-1}{\phi(p^a)}\r\rfloor+n-j$$
and therefore
$$a_j\gs\l\lfloor\f{n-j-1}{\phi(p^a)}\r\rfloor+n-j
+\l\lfloor\f{j-p^{a-1}-lp^a}{\phi(p^a)}\r\rfloor
-\l\lfloor\f{n-p^{a-1}-lp^a}{\phi(p^a)}\r\rfloor\gs n-j-1$$
by applying Lemma 2.1.

Combining the above we immediately obtain the desired result. \qed

\proclaim{Lemma 3.3} Let $p$ be a prime, $n\in\Z^+$, $r\in\Z$ and
$s,t\in\{0,\ldots,p-1\}$ with $s\not=p-1$. If $s<t$ then
$$\mm{pn+s}t{n-1,p}\eq(-1)^{n+s}\f n{t\bi{t-1}s}\ \ (\mo\ p).\tag3.3$$
If $s\gs t$, then
$$\mm{pn+s}t{n-1,p}\eq(-1)^{n+t}n\bi st\f{\sigma_{st}}p\ \ (\mo\ p),\tag3.4$$
where
$$\sigma_{st}=1+(-1)^p\f{\prod_{1\ls i\ls p,\, i\not=p-t}(p(n-1)+t+i)}
{\prod_{1\ls i\ls p,\, i\not=p-(s-t)}(s-t+i)}\eq1+(-1)^p\eq0\ (\mo\ p).\tag3.5$$
\endproclaim
\Proof. Clearly
$$\align \mm{pn+s}t{n-1,p}=&p^{-\l\lfloor\f{pn+s-1-(n-1)p}{p-1}\r\rfloor}
\sum_{k\eq t\,(\mo\ p)}(-1)^k\bi{pn+s}k\bi{(k-t)/p}{n-1}
\\=&\f{(-1)^{pn+t}}p\bi{pn+s}{pn+t}\bi n{n-1}
\\&+\f{(-1)^{p(n-1)+t}}p\bi{pn+s}{p(n-1)+t}\bi{n-1}{n-1}.
\endalign$$

{\tt Case 1}. $s<t$. In this case,  $d=t-1-s\gs0$ and
$$\align \mm{pn+s}t{n-1,p}=&\f{(-1)^{p(n-1)+t}}p
\prod_{i=0}^s\f{pn+i}{p(n-1)+t-i}\cdot\bi{p(n-1)+p-1}{p(n-1)+d}
\\=&\f{(-1)^{p(n-1)+t}n}{p(n-1)+t}
\prod_{i=1}^s\f{pn+i}{p(n-1)+t-i}\cdot\bi{p(n-1)+p-1}{p(n-1)+d}
\\\eq&(-1)^{n-1+t}\f{n\times s!}{\prod_{i=0}^s(t-i)}\bi{p-1}{d}
\ \ \ (\t{by Lucas' theorem})
\\\eq&(-1)^{n-s}\f{n}{t\bi{t-1}s}\ \ (\mo\ p).
\endalign$$

{\tt Case 2}. $s\gs t$.
Note that $$\sigma_{st}\eq1+(-1)^p\f{(p-1)!}{(p-1)!}\eq1+(-1)^p\eq0\ (\mo\ p)$$
and
$$\align \mm{pn+s}t{n-1,p}=&\f{(-1)^{pn+t}}p
\bi{pn+s}{pn+t}\(n+(-1)^p\prod_{i=1}^p\f{p(n-1)+t+i}{s-t+i}\)
\\=&(-1)^{pn+t}\f np\bi{pn+s}{pn+t}\sigma_{st}.
\endalign$$
Therefore
$$\mm{pn+s}t{n-1,p}\eq(-1)^{n+t}n\bi st\f{\sigma_{st}}p\quad (\mo\ p)$$
by Lucas' theorem.

The proof of Lemma 3.3 is now complete. \qed

\medskip
\noindent{\it Proof of Theorem 1.1}. If $n=0$ or $s=p-1$ or
$\phi(p^a)\nmid n-(l+1)p^{a-1}$, then (1.1) holds by Lemma 3.2.

Now we suppose that $n>0$, $s\not=p-1$ and $\phi(p^a)\mid
n-(l+1)p^{a-1}$. Then $p^{a-1}\mid n$, and hence $p\mid n$ since
$a\gs2$. Therefore $\mm{pn+s}t{n-1,p}\eq0\ (\mo\ p)$ by Lemma 3.3,
and thus we have (1.1) by (3.2).

This concludes the proof. \qed

\heading{4. Proofs of Theorem 1.2 and Corollary 1.3}\endheading

\proclaim{Lemma 4.1} Let $p$ be a prime, and let $a\in\Z^+$, $l\in\N$ and $r\in\Z$.
Then, for any $n\in\N$ with $n\eq l\ (\mo\ p-1)$, we have
$$\mm nr{l,p}\eq\cases(-1)^{\f{n-l}{p-1}-1}\bi{\f{n-l}{p-1}-1}l\ (\mo\ p)&\t{if}\ n>l,
\\0\ (\mo\ p)&\t{if}\ n\ls l.\endcases\tag4.1$$
\endproclaim
\Proof. We use induction on $m=(n-l)/(p-1)$.

If $m\ls l$ (i.e., $n\ls lp$), then $\lfloor(n-lp-1)/(p-1)\rfloor<0$,
and hence
$$\mm nr{l,p}=p^{-\lfloor\f{n-lp-1}{p-1}\rfloor}
\sum_{k\eq r\,(\mo\ p)}(-1)^k\bi nk\bi{(k-r)/p}l
\eq0\ (\mo\ p)$$
which yields (4.1).
If $l<m\ls 1$, then $l=0$ and $m=1$, hence $n=p-1$ and
$$\mm nr{l,p}\eq\sum_{k\eq r\,(\mo\ p)}\bi{p-1}k(-1)^k
\eq 1=(-1)^{m-1}\bi{m-1}l\ (\mo\ p).$$
Thus the desired result always holds in the case $m\ls\max\{l,1\}$.

Now let $m>\max\{l,1\}$ and
assume that whenever $l_*,n_*\in\N$ and $(n_*-l_*)/(p-1)=m-1>0$ we have
$$\mm{n_*}i{l_*,p}=(-1)^{\f{n_*-l_*}{p-1}-1}\bi{\f{n_*-l_*}{p-1}-1}{l_*}
=(-1)^{m}\bi{m-2}{l_*}\ (\mo\ p)$$
for all $i\in\Z$.

For $n'=n-(p-1)$ clearly $(n'-l)/(p-1)=m-1\gs \max\{l,1\}$.
By the induction hypothesis,
$\mm{n'}{i}{l,p}\eq(-1)^{m}\bi {m-2}l\ (\mo\ p)$
for each $i\in\Z$.
In view of the Chu-Vandermonde convolution identity (cf. [GKP, (5.27)]),
$$\bi nk=\sum_{j=0}^{p-1}\bi{p-1}j\bi{n'}{k-j}$$
for every $k=0,1,2,\ldots$. Therefore
$$\align\mm nr{l,p}=&p^{-\l\lfloor\f{n-lp-1}{p-1}\r\rfloor}
\sum_{j=0}^{p-1}\bi{p-1}j\sum_{p\mid k-r}(-1)^k\bi{n'}{k-j}\bi{(k-r)/p}l
\\=&\sum_{j=0}^{p-1}\bi{p-1}j\f{(-1)^j}p\mm{n'}{r-j}{l,p}
\\=&\sum_{j=0}^{p-1}\bi{p-1}j(-1)^j\f{\mm{n'}{r-j}{l,p}-(-1)^{m}\bi {m-2}l}p,
\endalign$$
since $\sum_{j=0}^{p-1}\bi{p-1}j(-1)^j=(1-1)^{p-1}=0$. Thus
$$\mm nr{l,p}\eq\sum_{j=0}^{p-1}\f{\mm{n'}{r-j}{l,p}
-(-1)^{m}\bi {m-2}l}p\ \ (\mo\ p).$$

 Observe that
$$\align &p^{\lfloor\f{n'-lp-1}{p-1}\rfloor}
\sum_{j=0}^{p-1}\mm{n'}{r-j}{l,p}
\\=&\sum_{j=0}^{p-1}\sum_{k\eq r-j\,(\mo\ p)}(-1)^k
\bi{n'}{k}\bi{(k-(r-j))/p}l
\\=&\sum_{k=0}^{n'}(-1)^k\bi{n'}{k}\bi{\lfloor(k-r+p-1)/p\rfloor}l
\\=&\cases\sum_{k=0}^{n'}(-1)^k\bi{n'}k=(1-1)^{n'}=0&\t{if}\ l=0,
\\-\sum_{k\eq r\,(\mo\ p)}(-1)^k\bi{n'-1}{k}\bi{(k-r)/p}{l-1}&\t{if}\ l>0,
\endcases
\endalign$$
where we have applied Lemma 2.1 of Sun [S06] to get the last equality.
Also,
$$\l\lfloor\f{n'-1-(l-1)p-1}{p-1}\r\rfloor=\l\lfloor\f{n'-lp-1}{p-1}\r\rfloor+1$$
and
$$ \f{n'-1-(l-1)}{p-1}=m-1.$$
Therefore
$$\align&\f1p\sum_{j=0}^{p-1}\mm{n'}{r-j}{l,p}
=\cases0&\t{if}\ l=0,\\-\mm{n'-1}{r}{l-1,p}&\t{if}\ l>0,\endcases
\\\eq&(-1)^{m-1}
\bi{m-2}{l-1}\ (\mo\ p)\ \ \t{(by the induction hypothesis)}.
\endalign$$

Combining the above we finally obtain that
$$\align \mm nr{l,p}\eq&\f1p\sum_{j=0}^{p-1}\mm{n'}{r-j}{l,p}-(-1)^{m}\bi {m-2}l
\\\eq&(-1)^{m-1}\bi {m-2}{l-1}+(-1)^{m-1}\bi {m-2}l
\\\eq&(-1)^{m-1}\bi{m-1}l\ \ (\mo\ p).
\endalign$$
This concludes the induction proof. \qed

\medskip
\noindent{\it Proof of Theorem 1.2}. By Lemma 3.2, if $s=p-1$, or
$\phi(p)=p-1$ does not divide $n-l-1$, then (1.2) holds. If
$s\not=p-1$ and $p\mid n$, then we also have (1.2) by Lemmas 3.2
and 3.3. Below we assume that $s\not=p-1$, $p-1\mid n-l-1$ and
$p\nmid n$.

When $s=2t$, clearly
$$\align\sigma_{st}=&1+(-1)^p\prod\Sb 1\ls i\ls p\\i\not=p-t\endSb\l(1+\f{p(n-1)}{t+i}\r)
\\\eq&1+(-1)^p\(1+p(n-1)\sum\Sb 1\ls i\ls p\\i\not=p-t\endSb\f1{t+i}\)\eq p\da_{p,2}\ \ (\mo\ p^2),
\endalign$$
for, $n$ is odd if $p=2$, and
$$\sum\Sb 1\ls i\ls p\\i\not= p-t\endSb\f1{t+i}\eq\sum_{k=1}^{p-1}\f1k
=\sum_{k=1}^{(p-1)/2}\l(\f1k+\f1{p-k}\r)\eq0\ \ (\mo\ p)$$ if
$p\not=2$. Therefore, in the case $s=2t$ and $p\not=2$, we have
(1.2) by Lemmas 3.2 and 3.3.

Now we consider the case $s<t$.
By Lemmas 3.2, 3.3 and 4.1,
$$\align\mm{pn+s}{pr+t}{l,p^2}\eq&(-1)^{n-1}\mm{pn+s}{t}{n-1,p}\mm{n-1}r{l,p}
\\\eq&(-1)^{n-1}(-1)^{n+s}\f n{t\bi{t-1}s}
\\&\times\cases(-1)^{\f{(n-1)-l}{p-1}-1}\bi{\f{n-1-l}{p-1}-1}l\ (\mo\ p)&\t{if}\ n-1>l,
\\0\ (\mo\ p)&\t{if}\ n-1\ls l.\endcases
\endalign$$

In view of the above we have completed the proof of Theorem 1.2. \qed

\medskip
\noindent{\it Proof of Corollary 1.3}.
We just modify the third case in the proof of Theorem 1.5 of [SD].
The only thing we require is that in the case
$n>0$ and $n\eq r\eq0\ (\mo\ p)$ we still have
$$T_{0,2}^{(p)}(n,r)=\f{p^{\lfloor\f{n/p-1}{p-1}\rfloor}}{(n/p)!}\mm nr{0,p^2}
\eq \f{p^{\lfloor\f{n_0-1}{p-1}\rfloor}}{n_0!}\mm {n_0}{r_0}{0,p}
=T_{0,1}^{(p)}(n_0,r_0)\ (\mo\ p)$$
where $n_0=n/p$ and $r_0=r/p$. Note that
$$\ord_p(n_0!)=\sum_{i=1}^{\infty}\l\lfloor\f{n_0}{p^i}\r\rfloor
<\sum_{i=1}^{\infty}\f{n_0}{p^i}=\f{n_0}{p-1}$$
and thus $\ord_p(n_0!)\ls \lfloor(n_0-1)/(p-1)\rfloor$.

If $p\not=2$, then by applying (1.2) with $l=s=t=0$ we find that
$$\mm nr{0,p^2}=\mm{pn_0}{pr_0}{0,p^2}\eq\mm {n_0}{r_0}{0,p}\ (\mo\ p)$$
and so $T_{0,2}^{(p)}(n,r)\eq T_{0,1}^{(p)}(n_0,r_0)\ (\mo\ p)$.
The last congruence also holds when $p=2$, because by Lemma 4.2 of [SD] we have
$$2\nmid T_{0,2}^{(2)}(n,r)
\iff n=2n_0\ \t{is a power of}\ 2
\iff 2\nmid T_{0,1}^{(2)}(n_0,r_0).$$
This concludes the proof. \qed

\heading{5. Proofs of Theorems 1.4 and 1.5}\endheading

\noindent{\it Proof of Theorem 1.4}.  By Lemma 3.2 of [SD] and its proof, if $j\in\N$ then
$$\sum_{k\eq0\,(\mo\ p)}(-1)^k\bi{pn}k\bi{k/p}j
=\sum_{j\ls k\ls n}(-1)^{pk}\bi{pn}{pk}\bi kj$$
is congruent to
$$\sum_{j\ls k\ls n}(-1)^k\bi nk\bi kj
=\bi nj\sum_{k\gs j}(-1)^k\bi{n-j}{k-j}=(-1)^j\,\da_{j,n}$$
modulo $p^{2\ord_p(n)+1+\da}$.
Therefore
$$\ord_p\(\mm{pn}0{j,p}\)\gs2\ord_p(n)+1+\da-\l\lfloor\f{pn-jp-1}{p-1}\r\rfloor$$
for any $j\in\N$ with $j\not=n$.
As in the proof of Lemma 3.2,
$$\mm{pn}{pr}{l,p^{a+1}}=(-1)^n(-1)^{pn}\mm nr{l,p^a}
+\sum_{0\ls j<n}(-1)^jp^{a_j}\mm jr{l,p^a}\mm{pn}0{j,p}$$
where $a_j\in\Z$ and $a_j\gs n-j-1$.

Let $m$ be the least integer greater than or equal to $\f{p-1}p(2\ord_p(n)+\da)$.
Then $m-1<\f{p-1}p(2\ord_p(n)+\da)$ and hence
$$m+\l\lfloor\f{m-1}{p-1}\r\rfloor=\l\lfloor\f{p(m-1)}{p-1}\r\rfloor+1\ls 2\ord_p(n)+\da.$$
For $0\ls j<n$, if $n-j\gs m+1$ then $a_j\gs n-j-1\gs m$; if $n-j\ls m$
then
$$\align a_j+\ord_p\(\mm{pn}0{j,p}\)
\gs&n-j-1+2\ord_p(n)+1+\da-\l\lfloor\f{p(n-j)-1}{p-1}\r\rfloor
\\=&2\ord_p(n)+\da-\l\lfloor\f{n-j-1}{p-1}\r\rfloor
\\\gs&2\ord_p(n)+\da-\l\lfloor\f{m-1}{p-1}\r\rfloor\gs m.
\endalign$$

Combining the above we get that
$$\ord_p\(\mm{pn}{pr}{l,p^{a+1}}-(-1)^{(p-1)n}\mm nr{l,p^a}\)
\gs m\gs\f{p-1}p(2\ord_p(n)+\da).$$
If $(p-1)n$ is odd, then $p=2$ and $2\nmid n$, hence $2\ord_p(n)+\da=0$.
So (1.5) holds. \qed

\medskip
\noindent{\it Proof of Theorem 1.5}.
We use induction on $a$.

When $a=1$, the desired result follows from Lemma 4.1.

In the case $a=2$, by Theorem 1.2 and Lemma 4.1, we have
$$\align\mm nr{l,p^2}=&\mm {p(l+m(p-1))+p-1}{p\lfloor r/p\rfloor+\{r\}_p}{l,p^2}
\\\eq&(-1)^{\{r\}_p}\bi{p-1}{\{r\}_p}\mm{l+m(p-1)}{\lfloor r/p\rfloor}{l,p}
\eq\mm{l+m(p-1)}{\lfloor r/p\rfloor}{l,p}
\\\eq&(-1)^{m-1}\bi{m-1}l\ \ (\mo\ p).
\endalign$$

Now let $a>2$ and assume Theorem 1.5 with $a$ replaced by $a-1$. Then,
with helps of Theorem 1.1 and the induction hypothesis, we have
$$\align\mm nr{l,p^a}=&\mm{p^{a-1}(l+m(p-1)+1)-1}r{l,p^a}
\\=&\mm{p(p^{a-2}(l+m(p-1)+1)-1)+(p-1)}{p\lfloor r/p\rfloor+\{r\}_p}{l,p^a}
\\\eq&(-1)^{\{r\}_p}\bi{p-1}{\{r\}_p}
\mm{p^{a-2}(l+m(p-1)+1)-1}{\lfloor r/p\rfloor}{l,p^{a-1}}
\\\eq&\mm{(l+1)p^{a-2}-1+m\phi(p^{a-1})}{\lfloor r/p\rfloor}{l,p^{a-1}}
\\\eq&(-1)^{m-1}\bi{m-1}l\ \ (\mo\ p).
\endalign$$
This concludes the induction step and we are done. \qed

\enddocument